February 11, 2018# GREAT CIRCLE FIBRATIONS and CONTACT STRUCTURES on the 3-SPHERE

## Herman Gluck

Given any smooth fibration of the unit 3-sphere by great circles, we show that the distribution of 2-planes orthogonal to the great circle fibres is a tight contact structure, a fact well known in the special case of the Hopf fibrations. The proof expresses hypothesis and conclusion as differential inequalities involving functions on disks transverse to the fibres, and shows that one inequality implies the other.

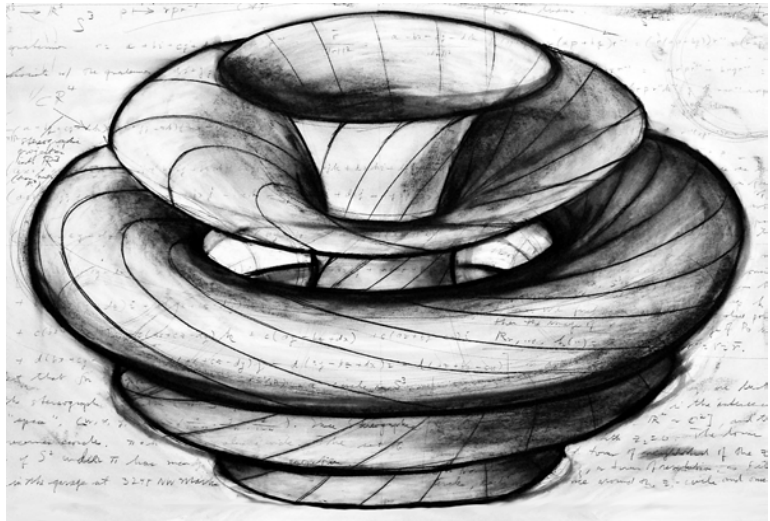

Hopf fibration of 3-sphere by great circles
Lun-Yi Tsai  Charcoal and graphite on paper  2006

**MAIN THEOREM.** Let $F$ be a smooth fibration of $S^3$ by great circles, and $\xi_F$ the distribution of tangent 2-planes orthogonal to the fibres of $F$. Let $A$ be a unit vector field tangent to the fibres of $F$, and $\alpha$ the dual one-form defined by $\alpha(V) = <A, V>$, so that $\xi_F = \ker \alpha$. Then $\alpha \wedge d\alpha$ is never 0, and therefore $\xi_F$ is a contact structure on $S^3$, which we show is tight.

**An overview of the proofs.**

**Overview of proof of the Main Theorem.** Start with the smooth fibration $F$ of $S^3$ by great circles, the orthogonal 2-plane distribution $\xi_F$ and the related one-form $\alpha$.

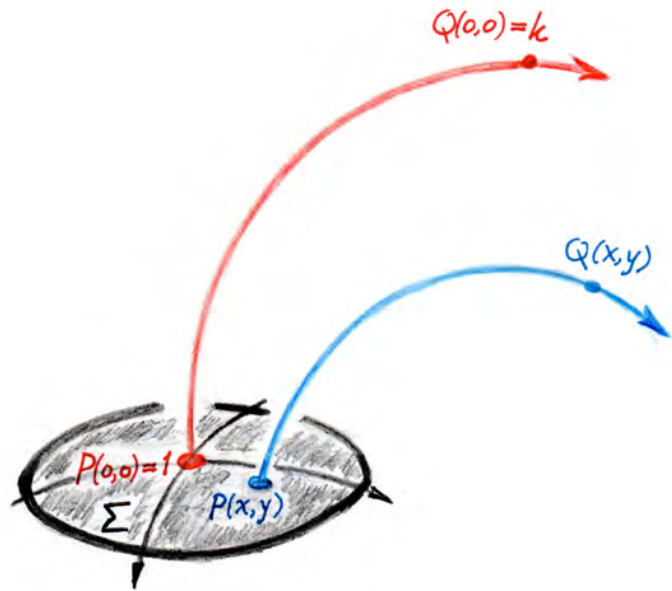

Regard $S^3$ as the space of unit quaternions, and focus on a random fibre of $F$, positioned to run through $1$ and $k$ and shown above in red. Consider a small round disk $\Sigma$ which lies on a great 2-sphere meeting this fibre orthogonally at the center of the disk.

At each point $P(x, y)$ on $\Sigma$, consider the fibre of $F$ shown in blue which pierces $\Sigma$ at that point. Follow this fibre for a distance $\pi/2$ to arrive at the point

$$Q(x, y) = P(x, y) \left( f(x, y)\, \mathbf{i} + g(x, y)\, \mathbf{j} + h(x, y)\, \mathbf{k} \right),$$

where $f^2 + g^2 + h^2 = 1$, and where we use quaternion multiplication. For the Hopf fibration $H$ we have $f \equiv g \equiv 0$ and $h \equiv 1$, so the functions $f$ and $g$ measure the deviation of $F$ from $H$.

**Proposition 1.** The constraint $(f_x - g_y)^2 < 4(1 + f_y)(1 - g_x)$ at the origin $(x, y) = (0, 0)$ is necessary and sufficient for the family $F$ of great circles on $S^3$ corresponding to the functions $f$, $g$ and $h$ to be a smooth fibration in some small tubular neighborhood of the given great circle fibre through $1$ and $k$.

**Proposition 2.** The constraint $(1 + f_y) + (1 - g_x) > 0$ at the origin $(x, y) = (0, 0)$ is necessary and sufficient for $\alpha \wedge d\alpha \neq 0$ there, and hence for $\xi_F$ to be a contact structure in a neighborhood of this point.

We then say why the first constraint implies the second, completing the proof of the main theorem.

## Acknowledgments.

Many thanks to Patricia Cahn, Haggai Nuchi and Jingye Yang for their thorough proof reading of the arguments here, and to Yasha Eliashberg and John Etnyre for their very thoughtful comments.

## Tools to be used in the proofs.

### (1) The Grassmann manifold $G_2R^4$.

An oriented great circle on the unit 3-sphere $S^3 \subset R^4$ spans an oriented 2-plane through the origin in $R^4$, and the set of these is the Grassmann manifold $G_2R^4$.

There are many ways to identify $G_2R^4$ with $S^2 \times S^2$; here is one which works well with our proofs. View $R^4$ as the space of quaternions, $S^3$ as the subspace of unit quaternions, and $S^2$ as its subspace of unit, purely imaginary quaternions. Then given an oriented 2-plane through the origin in $R^4$, and an ordered orthonormal basis $P$ and $Q$ for this 2-plane, we consider the pair $(Q\underline{P}, \underline{P}Q)$ in $S^2 \times S^2$, where the under-score indicates quaternionic conjugation.

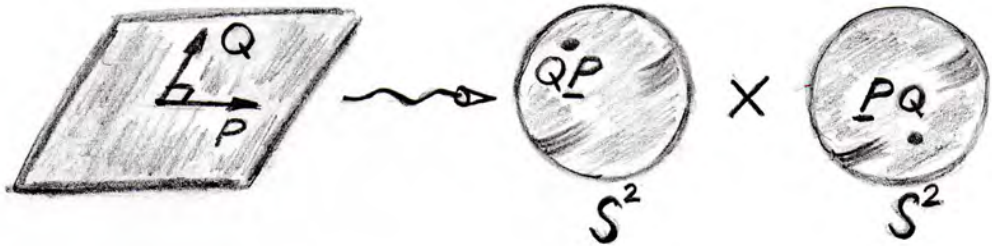

We saw in [DeTurck, Gluck, et al, 2013] that this provides a picture of $G_2R^4$ as $S^2 \times S^2$, isometric up to scale.

### (2) The screw sense of a fibration of $S^3$ by great circles.

The family of Hopf fibrations of $S^3$ by great circles consists of two disjoint copies of $S^2$, distinguished by the so-called "screw sense" of the fibres. A single Hopf fibration $H$ can be determined by the choice of a unit, purely imaginary quaternion, say $q = a\mathbf{i} + b\mathbf{j} + c\mathbf{k}$, and by the feature that the fibres of $H$ are invariant under multiplication by $q$. If we multiply by $q$ on the right, we call $H$ a **right Hopf fibration**, and if we multiply by $q$ on the left, then a **left Hopf fibration**. All the right Hopf fibrations lie in one copy of $S^2$, all the left ones in the other copy.

**(3) The space of all great circle fibrations of $S^3$ deformation retracts to its subspace of Hopf fibrations.**

We proved this in [Gluck - Warner, 1983].

In this way, all great circle fibrations of $S^3$ inherit a screw sense from the Hopf fibrations.

**(4) The base space of a great circle fibration of $S^3$ viewed in $G_2 R^4$.**

Let F be a fibration of $S^3$ by oriented great circles. Each fibre spans an oriented 2-plane through the origin in $R^4$, which corresponds to a point in the Grassmann manifold $G_2 R^4$. In this way, the base space $M_F$ of the fibration F appears as a submanifold $M_F$ of $G_2 R^4$. When we view $G_2 R^4$ as $S^2 \times S^2$, as above, the base space $M_F$ appears as the graph of a strictly distance-decreasing map $\varphi$ from one $S^2$ factor to the other.

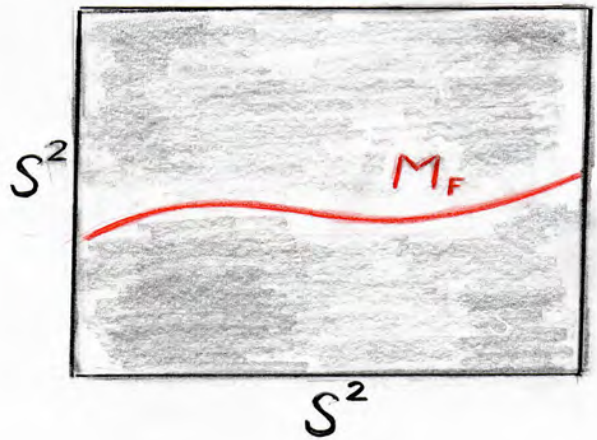

If the fibration F is smooth, then the norm $|d\varphi| < 1$, meaning that the linear map $d\varphi$ strictly decreases the lengths of all non-zero vectors. See [Gluck - Warner, 1983] for details.

## (5) Deforming a great circle fibration of $S^3$ to a Hopf fibration while fixing a selected fibre.

Since the space of all great circle fibrations of $S^3$ deformation retracts to its subspace of Hopf fibrations, we can start with any such fibration $F$ and then find a path of great circle fibrations $F_t$, $0 \le t \le 1$, such that $F_0 = F$ and $F_1 = $ some Hopf fibration $H$.

Near the end of the proofs of Propositions 1 and 2, we will start with such a great circle fibration $F$ and within it, a distinguished fibre, and want to choose the above path $F_t$ of great circle fibrations so that each one contains that distinguished fibre.

To do so, we enter into the proof of the above deformation retraction in [Gluck - Warner, 1983], and modify it as follows.

Suppose the initial fibration $F$ corresponds to the distance-decreasing map $\varphi: S^2 \to S^2$. In the proof, we saw that the image $\varphi(S^2)$ is a closed set in an open hemisphere of the range $S^2$. We noted there that $\varphi(S^2)$ must be contained in a unique smallest closed round disk in that hemisphere.

In the proof of deformation retraction, we shrunk that disk to its center, thereby deforming $\varphi$ to a constant map of $S^2$ to that center point, knowing that the graph of such a constant map is the base space of some Hopf fibration. We chose the center of the disk as the destination to guarantee uniformity of our process.

But now suppose that it is just a single great circle fibration $F$ that we are deforming, and that we have distinguished one of its fibres that we want to preserve during the deformation. Suppose that fibre corresponds to the point $(x, y) \in S^2 \times S^2$.

Then the point $y$ will lie inside the smallest closed round disk containing the image $\varphi(S^2)$, and we simply deform that disk to the point $y$ instead of to its center. In that way, all the fibrations $F_t$ in the deformation will contain the distinguished fibre corresponding to this point $(x, y)$, as desired.

## Set up for proving the Main Theorem.

We start with a smooth fibration F of $S^3$ by great circles, and let $\xi_F$ be the distribution of tangent 2-planes orthogonal to the fibres of F.

Since $S^3$ is simply connected, we can orient all these great circle fibres consistently in one direction or the other, and just choose one. Then we let A be a unit vector field tangent to the fibres of F, pointing in the direction of their chosen orientation.

Let $\alpha$ be the dual one-form defined by $\alpha(V) = <A, V>$ for any vector field V on $S^3$, so that the 2-planes $\xi_F$ are the kernels of $\alpha$ at each point.

We will show that $\alpha \wedge d\alpha$ is never 0, which tells us that $\xi_F$ is a contact structure on $S^3$.

Afterwards, we will explain why this contact structure is tight.

*****

We view $R^4$ as the space of quaternions $w + x\mathbf{i} + y\mathbf{j} + z\mathbf{k}$, and then view $S^3$ as the subspace of unit quaternions.

To show that $\alpha \wedge d\alpha$ is never zero, we pick at random a point of $S^3$ and will show that $\alpha \wedge d\alpha \neq 0$ there.

We turn the 3-sphere around so that this point is located at 1 and so that the fibre of F through 1 also goes through $\mathbf{k}$, and is oriented from 1 towards $\mathbf{k}$.

We let $\Sigma$ denote a small round disk lying on the great 2 sphere $S^3 \cap \{k = 0\}$ and centered at 1, which is orthogonal to the great circle fibre of F through this point.

We denote a typical point on $\Sigma$ by

$$P(x, y) = \left((1 - x^2 - y^2)^{1/2}, x, y, 0\right), \text{ with } x^2 + y^2 < \varepsilon.$$

Then $P(0, 0) = 1 = (1, 0, 0, 0)$ is the center of our disk $\Sigma$.

We want to name the fibre of F which passes through the point P(x, y) on $\Sigma$.

So we go out a distance $\pi/2$ along this fibre, in the direction of its orientation, and come to the point Q(x, y), which as a vector in $R^4$ is orthogonal to P(x, y).

Thus we can write $Q(x, y) = P(x, y) R(x, y)$, where R(x, y) is a unit purely imaginary quaternion, and where juxtaposition of P and R indicates quaternion multiplication.

We will assume that our fibration F has a right-handed screw sense.

We write

$$R(x, y) = f(x, y)\,\mathbf{i} + g(x, y)\,\mathbf{j} + h(x, y)\,\mathbf{k},$$

with $f^2 + g^2 + h^2 = 1$, and with $R(0, 0) = \mathbf{k}$, and hence $f(0, 0) = 0$, $g(0, 0) = 0$ and $h(0, 0) = 1$.

Then

$$Q(x, y) = P(x, y)\,R(x, y) = \left((1 - x^2 - y^2)^{1/2},\, x,\, y,\, 0\right)(0,\, f,\, g,\, h)$$

$$= \left(-xf - yg,\, (1 - x^2 - y^2)^{1/2} f + yh,\, (1 - x^2 - y^2)^{1/2} g - xh,\, (1 - x^2 - y^2)^{1/2} h + xg - yf\right)$$

is the point on the fibre of F through $P(x, y)$, a quarter circle ahead of it.

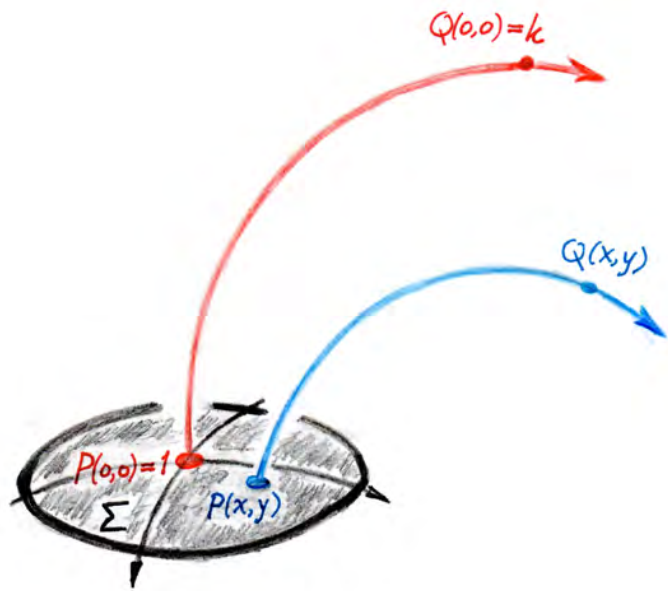

We think of the small round spherical disk $\Sigma$ as a "firing ground", from which points are fired that follow great circle trajectories that are fibres of F. The trajectory from the center $P(0, 0) = 1$ of the disk goes straight "up", orthogonal to the disk $\Sigma$, but the other trajectories take off at various angles to $\Sigma$.

The "firing solution" is given by a choice of the three functions $f(x, y)$, $g(x, y)$ and $h(x, y)$, which must be constrained so that the great circle trajectories, that is, the fibres of F, do not collide with one another.

Figuring out this constraint is the task of Proposition 1.

**Proof of Proposition 1.**

The points P(x, y) and Q(x, y) are $\pi/2$ apart on the blue fibre of F shown in the figure, and hence provide an ordered orthonormal basis for the oriented 2-plane through the origin in $R^4$ spanned by that fibre.

As mentioned in the Tools section, this oriented 2-plane corresponds to the pair (Q P , P Q) in the $S^2 \times S^2$ picture of the Grassmann manifold $G_2 R^4$.

Then the set

$$\Sigma_F = \{ (Q(x, y) \underline{P(x, y)}, \underline{P(x, y)} Q(x, y)) : x^2 + y^2 < \varepsilon \}$$

is the base space for the thin tube of fibres of F passing through the small disk $\Sigma$.

The constraint on this "firing solution" comes from the fact that $\Sigma_F \subset S^2 \times S^2$ must be the graph of a strictly distance-decreasing smooth map from an open set in one $S^2$ factor to the other $S^2$ factor, as mentioned in the Tools section.

We compute that

$$Q(x, y) \underline{P(x, y)} = \big( 0, \, 2xyg + (1 - 2y^2)f + 2(1 - x^2 - y^2)^{1/2} yh,$$

$$2xyf + (1 - 2x^2)g - 2(1 - x^2 - y^2)^{1/2} xh,$$

$$(1 - 2x^2 - 2y^2)h + 2(1 - x^2 - y^2)^{1/2} xg - 2(1 - x^2 - y^2)^{1/2} yf \big),$$

while

$$\underline{P(x, y)} Q(x, y) = R(x, y) = (0, f, g, h).$$

We let $D_\varepsilon$ denote the open disk $D_\varepsilon = \{(x, y): x^2 + y^2 < \varepsilon\}$ in the xy-plane, and consider the maps $\pi_1: D_\varepsilon \to S^2$ and $\pi_2: D_\varepsilon \to S^2$ given by

$$\pi_1(x, y) = Q(x, y) \underline{P(x, y)} \text{ and } \pi_2(x, y) = \underline{P(x, y)} Q(x, y).$$

We must compute the differential of $\pi_2 \pi_1^{-1}$ at the origin $x = 0$, $y = 0$, and discover the constraint on the choice of $f(x, y)$, $g(x, y)$ and $h(x, y)$ which guarantees that this linear map is strictly distance-decreasing.

To this end, we compute that at the origin $(x, y) = (0, 0)$, we have

$$\partial\pi_1/\partial x = (f_x, g_x - 2h, h_x + 2g) \text{ and } \partial\pi_1/\partial y = (f_y + 2h, g_y, h_y - 2f).$$

Recall that $f^2 + g^2 + h^2 = 1$, and then differentiate this with respect to $x$ to get

$$f f_x + g g_x + h h_x = 0.$$

At the origin we have $f = 0$, $g = 0$ and $h = 1$, and so conclude that there we also have $h_x = 0$. And likewise, we have $h_y = 0$ there.

Inserting this information above, we see that at the origin we have

$$\partial\pi_1/\partial x = (f_x, g_x - 2, 0) \text{ and } \partial\pi_1/\partial y = (f_y + 2, g_y, 0).$$

Thus at the origin, we have that $d\pi_1$ is given by the 2 x 2 matrix

$$\begin{matrix} f_x & f_y + 2 \\ g_x - 2 & g_y \end{matrix} \quad,$$

where we have dropped the last row of zeros.

And likewise, but more easily, we get that at the origin, $d\pi_2$ is given by the 2 x 2 matrix

$$\begin{matrix} f_x & f_y \\ g_x & g_y \end{matrix} \quad.$$

For simplicity, we introduce the abbreviations

$$a = f_x, \; b = f_y + 2, \; c = g_x - 2 \text{ and } d = g_y,$$

and then put $\Delta = ad - bc$.

With these abbreviations, we compute that at the origin, $d(\pi_2 \pi_1^{-1}) = (d\pi_2)(d\pi_1)^{-1}$ is given by the 2 x 2 matrix $M = 1/\Delta$ times

$$\begin{matrix} \Delta + 2c & -2a \\ 2d & \Delta - 2b. \end{matrix}$$

The question now is, what are the constraints on the functions f, g and h which guarantee that M is strictly distance-decreasing?

The condition that the linear map represented by the 2 x 2 matrix M be strictly distance-decreasing is equivalent to the constraint that the largest eigenvalue of $M^T M$ be < 1, see [Golub - Van Loan, 1996], and this in turn translates to the constraint that

$$|M|^2 < 1 + (\det M)^2,$$

where $|M|^2$ is the sum of the squares of the entries of M.

Now the entries of M were given above in terms of $a = f_x$, $b = f_y + 2$, $c = g_x - 2$ and $d = g_y$, and we compute that the above constraint translates to

$$(f_x - g_y)^2 < 4(1 + f_y)(1 - g_x)$$

at the origin $(x, y) = (0, 0)$.

Thus we have proved

**Proposition 1.** The constraint $(f_x - g_y)^2 < 4(1 + f_y)(1 - g_x)$ at the origin $(x, y) = (0, 0)$ is necessary and sufficient for the family F of great circles on $S^3$ corresponding to the "firing solution" (f, g, h) to be a smooth fibration in some small tubular neighborhood of the given great circle fibre through 1 and k.

**Remark.** In the constraint inequality above, the two factors $1 + f_y$ and $1 - g_x$ could either both be positive or both be negative. We claim they both must be positive.

To see this, consider the Hopf fibration where the great circle through any point $P$ on $S^3$ also goes through the quaternionic product $P\,k$, just like our base circle through $P = 1$.

The firing data for that Hopf fibration is given by $f \equiv 0$, $g \equiv 0$, $h \equiv 1$, and for that data we have

$$(1 + f_y) = 1 > 0 \quad \text{and} \quad (1 - g_x) = 1 > 0$$

at the origin $(x, y) = (0, 0)$.

Now this Hopf fibration, just like our fibration $F$, has a right-handed screw sense. Hence, as seen in the Tools section, we can connect these two fibrations by a one-parameter family of great circle fibrations, all of which contain the same great circle through $1$ and $k$.

Then the inequality $(f_x - g_y)^2 < 4\,(1 + f_y)\,(1 - g_x)$ of Proposition 1 above implies that

$$(1 + f_y)\,(1 - g_x) > 0$$

for each of the great circle fibrations in our one-parameter family, and since the two factors $(1 + f_y)$ and $(1 - g_x)$ each start out positive at the Hopf fibration, neither one of them can vanish during the transition to our fibration, and so we must have both

$$1 + f_y > 0 \quad \text{and} \quad 1 - g_y > 0$$

throughout the transition, and in particular for our given great circle fibration $F$.

## Proof of Proposition 2.

Once again, we start with a smooth fibration $F$ of $S^3$ by great circles, and let $\xi_F$ be the distribution of tangent 2-planes orthogonal to the fibres of $F$. We let $A$ be a unit vector field tangent to the fibres of $F$, and $\alpha$ the dual one-form defined by $\alpha(V) = <A, V>$, so that $\xi_F = \ker \alpha$. We must show that $\alpha \wedge d\alpha$ is never zero, and that therefore $\xi_F$ is a contact structure on $S^3$.

Now corresponding to the firing solution $(f, g, h)$ we get at least a fibration of a small neighborhood of $P(0, 0) = 1$ by arcs of great circles, and so now seek the constraint on $f, g$ and $h$ which guarantees that $\alpha \wedge d\alpha$ is non-zero at this point, which corresponds to the origin $x = 0$ and $y = 0$ of our local coordinate system.

A small neighborhood of this point can be coordinatized by

$$S(x, y, t) = P(x, y) \cos t + Q(x, y) \sin t .$$

We now set out to prove

**Proposition 2.** The constraint $(1 + f_y) + (1 - g_x) > 0$ at the origin $(x, y) = (0, 0)$ is necessary and sufficient for $\alpha \wedge d\alpha \neq 0$ at that randomly chosen point, and hence for $\xi_F$ to be a contact structure in a neighborhood of the point.

To begin the argument, we recall that

$$P(x, y) = \left((1 - x^2 - y^2)^{1/2}, x, y, 0\right), \text{ with } x^2 + y^2 < \varepsilon ,$$

and that

$$Q(x, y) = \left(-xf - yg, (1 - x^2 - y^2)^{1/2} f + yh, (1 - x^2 - y^2)^{1/2} g - xh, (1 - x^2 - y^2)^{1/2} h + xg - yf\right) .$$

Then on our small coordinatized neighborhood of $P(0, 0) = 1$, the coordinate vector fields are given by

$$\partial/\partial x = \partial S/\partial x = P_x \cos t + Q_x \sin t$$

$$\partial/\partial y = \partial S/\partial y = P_y \cos t + Q_y \sin t$$

$$\partial/\partial t = \partial S/\partial t = -P \sin t + Q \cos t .$$

Within our coordinatized neighborhood, the unit vector field A tangent to the fibres of F is given by

$$A = \partial/\partial t = \partial S/\partial t = -P \sin t + Q \cos t.$$

The dual 1-form $\alpha$ to $A = \partial/\partial t$ must satisfy

$$\alpha(\partial/\partial t) = 1, \quad \alpha(\partial/\partial x) = \langle \partial/\partial t, \partial/\partial x \rangle, \quad \alpha(\partial/\partial y) = \langle \partial/\partial t, \partial/\partial y \rangle.$$

Let us temporarily write

$$\alpha = a\, dx + b\, dy + c\, dt,$$

where the one-forms dx, dy, dt are dual to the vector fields $\partial/\partial x$, $\partial/\partial y$, $\partial/\partial t$ in the sense that

$$dx(\partial/\partial x) = 1, \quad dx(\partial/\partial y) = 0, \quad \text{etc.}$$

Then

$$\alpha = \langle \partial/\partial t, \partial/\partial x \rangle\, dx + \langle \partial/\partial t, \partial/\partial y \rangle\, dy + \langle \partial/\partial t, \partial/\partial t \rangle\, dt.$$

We have

$$a = \langle \partial/\partial t, \partial/\partial x \rangle = \langle -P \sin t + Q \cos t, P_x \cos t + Q_x \sin t \rangle$$

$$= (-\langle P, P_x \rangle + \langle Q, Q_x \rangle) \sin t \cos t - \langle P, Q_x \rangle \sin^2 t + \langle Q, P_x \rangle \cos^2 t$$

$$= -\langle P, Q_x \rangle \sin^2 t + \langle Q, P_x \rangle \cos^2 t,$$

since $\langle P, P \rangle = \langle Q, Q \rangle = 1$ implies $\langle P, P_x \rangle = \langle Q, Q_x \rangle = 0$.

Likewise we have

$$b = \langle \partial/\partial t, \partial/\partial y \rangle = (-\langle P, P_y \rangle + \langle Q, Q_y \rangle) \sin t \cos t - \langle P, Q_y \rangle \sin^2 t + \langle Q, P_y \rangle \cos^2 t$$

$$= -\langle P, Q_y \rangle \sin^2 t + \langle Q, P_y \rangle \cos^2 t,$$

since $\langle P, P_y \rangle = \langle Q, Q_y \rangle = 0$. And of course we have

$$c = \langle \partial/\partial t, \partial/\partial t \rangle = 1.$$

These are the values of $a$, $b$ and $c$, which can be further broken down into combinations of $f$, $g$, $h$ and their derivatives ... but for simplicity, we continue to write

$$\alpha = a\, dx + b\, dy + dt.$$

Then

$$d\alpha = (b_x - a_y)\, dx \wedge dy + a_t\, dt \wedge dx + b_t\, dt \wedge dy,$$

and hence

$$\alpha \wedge d\alpha = (-a\, b_t + b\, a_t + b_x - a_y)\, dx \wedge dy \wedge dt.$$

We already know that

$$a = -\langle P, Q_x \rangle \sin^2 t + \langle Q, P_x \rangle \cos^2 t,$$

$$b = -\langle P, Q_y \rangle \sin^2 t + \langle Q, P_y \rangle \cos^2 t,$$

and then compute that

$$a_t = -2(\langle P, Q_x \rangle + \langle Q, P_x \rangle) \sin t \cos t$$

$$b_t = -2(\langle P, Q_y \rangle + \langle Q, P_y \rangle) \sin t \cos t$$

$$a_y = (-\langle P_y, P_x \rangle - \langle P, P_{xy} \rangle + \langle Q_y, Q_x \rangle + \langle Q, Q_{xy} \rangle) \sin t \cos t$$

$$\quad - (\langle P_y, Q_x \rangle + \langle P, Q_{xy} \rangle) \sin^2 t + (\langle Q_y, P_x \rangle + \langle Q, P_{xy} \rangle) \cos^2 t$$

$$b_x = (-\langle P_x, P_y \rangle - \langle P, P_{yx} \rangle) + \langle Q_x, Q_y \rangle + \langle Q, Q_{yx} \rangle) \sin t \cos t$$

$$\quad - (\langle P_x, Q_y \rangle + \langle P, Q_{yx} \rangle) \sin^2 t + (\langle Q_x, P_y \rangle + \langle Q, P_{yx} \rangle) \cos^2 t.$$

At $t = 0$, all the $\sin t$ terms drop out, and we get

$a = \langle Q, P_x \rangle$

$b = \langle Q, P_y \rangle$

$a_t = 0$

$b_t = 0$

$a_y = \langle Q_y, P_x \rangle + \langle Q, P_{xy} \rangle$

$b_x = \langle Q_x, P_y \rangle + \langle Q, P_{yx} \rangle.$

Then the coefficient of $dx \wedge dy \wedge dt$ in the expression for $\alpha \wedge d\alpha$ is given by

$$-a\, b_t + b\, a_t + b_x - a_y = b_x - a_y$$

$$= (<Q_x, P_y> + <Q, P_{yx}>) - (<Q_y, P_x> + <Q, P_{xy}>)$$

$$= <Q_x, P_y> - <Q_y, P_x>,$$

thanks to the cancellation of the two higher derivative terms.

Having already set $t = 0$, we now set $x = 0$ and $y = 0$ so that we may see the value of $\alpha \wedge d\alpha$ at the point $P(0, 0) = 1$.

At this point, we have

$$P = (1, 0, 0, 0) = 1$$

$$Q = (0, 0, 0, 1) = k$$

$$P_x = (0, 1, 0, 0) = i$$

$$P_y = (0, 0, 1, 0) = j$$

$$Q_x = (0, f_x, g_x - 1, 0) = f_x\, i + (g_x - 1)\, j$$

$$Q_y = (0, f_y + 1, g_y, 0) = (f_y + 1)\, i + g_y\, j,$$

with the last four lines above computed from the explicit formulas for $P(x, y)$ and $Q(x, y)$ given above, using the facts that at $P(0, 0)$ we have

$$f = 0,\ g = 0,\ h = 1,\ \text{and}\ h_x = 0,\ h_y = 0.$$

Inserting these values into the above expression for the coefficient $<Q_x, P_y> - <Q_y, P_x>$ of $dx \wedge dy \wedge dt$ in $\alpha \wedge d\alpha$, we get

$$<Q_x, P_y> - <Q_y, P_x> = (g_x - 1) - (f_y + 1) = -\big((1 + f_y) + (1 - g_x)\big),$$

and hence

$$\alpha \wedge d\alpha = -\big((1 + f_y) + (1 - g_x)\big)\, dx \wedge dy \wedge dt$$

at the point $P(0, 0) = 1$.

Thus the constraint

$$(1 + f_y) + (1 - g_x) \neq 0$$

at the origin $(x, y) = (0, 0)$ is necessary and sufficient for $\alpha \wedge d\alpha \neq 0$ there, and hence for $\xi_F$ to be a contact structure in a neighborhood of this point.

We can refine this constraint to read

$$(1 + f_y) + (1 - g_x) > 0$$

at $(x, y) = (0, 0)$ by an argument similar to that in the Remark following the proof of Proposition 1.

**Proof of the Main Theorem.**

We start with a smooth fibration $F$ of $S^3$ by great circles, and let $\xi_F$ be the distribution of tangent 2-planes orthogonal to the fibres of $F$. We let $A$ be a unit vector field tangent to the fibres of $F$, and $\alpha$ the dual one-form defined by $\alpha(V) = <A, V>$, so that $\xi_F = \ker \alpha$. We must show that $\alpha \wedge d\alpha$ is never zero, and that therefore $\xi_F$ is a contact structure on $S^3$.

So as before, we pick a point on $S^3$ at random, and will show that $\alpha \wedge d\alpha \neq 0$ there.

We view the 3-sphere $S^3$ as the space of unit quaternions, turn it around so that the given point is at the identity $1$ and so that the fibre of $F$ through $1$ also goes through $k$.

Then by Proposition 1, the firing solution $(f, g, h)$ for this fibration must satisfy the inequality

$$(f_x - g_y)^2 < 4(1 + f_y)(1 - g_x)$$

at our given point $P(0, 0) = 1$, and by the Remark following its proof, each of the factors $(1 + f_y)$ and $(1 - g_y)$ must be strictly positive.

But then we certainly have

$$(1 + f_y) + (1 - g_x) > 0,$$

and so by Proposition 2 we have that $\alpha \wedge d\alpha \neq 0$ at the given point.

This completes the proof of the Main Theorem, save that we must still explain why the contact sructure $\xi_F$ is tight.

## Proof that the contact structure $\xi_F$ is tight.

Let F be a smooth great circle fibration of $S^3$, and H the Hopf fibration to which it is connected in the deformation retraction of the space of all great circle fibrations of $S^3$ to its subspace of Hopf fibrations.

This deformation retraction provides a one-parameter family $F_t$ of such fibrations, which begins with F at $t = 0$ and ends with H at $t = 1$.

Then the corresponding contact forms $\alpha$ for F and $\alpha'$ for H can also be connected by a one-parameter family $\alpha_t$ of contact forms.

Hence by the Gray Stability Theorem [Gray, 1959; Geiges, 2008]], there is an isotopy $h_t$ of diffeomorphisms of $S^3$ with $h_0$ = identity and with $h_t{}^*(\alpha_0) = f(t)\,\alpha_t$, where f(t) is a real-valued function.

Thus the contact structures $\xi_F$ and $\xi_H$ are isotopic, meaning that there is a diffeomorphism h: $S^3 \to S^3$, isotopic to the identity, such that $dh(\xi_F) = \xi_H$.

Since $\xi_H$ is tight, so also is $\xi_F$ tight.

This completes the proof of our Main Theorem.

**Next steps.**

**(1)** As mentioned earlier, it is known that the space of great circle fibrations of $S^3$ deformation retracts to its subspace of Hopf fibrations.

It would be good to learn that the space of tight contact structures on $S^3$ deformation retracts to its subspace of standard ones (those orthogonal to Hopf fibrations). See [Eliashberg, 1992, 1993] for what is known.

In such a case, it would be sensible to ask if this deformation retraction to the subspace of standard contact structures on $S^3$ can be carried out so that at the same time, those which are orthogonal to great circle fibrations deform within themselves to the standard ones ... thus giving us a deformation retraction of pairs.

**(2)** Given a smooth fibration of $S^{2n+1}$ by great circles, is the distribution of 2n-planes orthogonal to the great circle fibres a contact structure?

Thanks to Peter McGrath for asking this question.

## References.

University of Pennsylvania
*gluck@math.upenn.edu*